\documentclass{article}

\usepackage[T1]{fontenc}
\usepackage[latin1]{inputenc}
\usepackage{lmodern}

\usepackage{graphicx}

\usepackage{amssymb}

\usepackage[intlimits]{amsmath}

\usepackage[all,v2]{xy}

\newtheorem{theorem}{Theorem}[section]

\newtheorem{lemma}[theorem]{Lemma}

\newtheorem{proposition}[theorem]{Proposition}

\newtheorem{definition}[theorem]{Definition\rm}

\newtheorem{remark}[theorem]{Remark}

\newtheorem{example}[theorem]{Example}

\newenvironment{proof}{\begin{trivlist}\item[]{\sc Proof.}}{\nolinebreak $\Box$ \end{trivlist}}

\usepackage{wasysym}
\usepackage{pifont}
\usepackage{stmaryrd}

\def\dar[#1]{\ar@<2pt>[#1]\ar@<-2pt>[#1]}

\DeclareMathOperator{\aut}{\mathfrak{aut}}

\DeclareMathOperator{\naive}{naive}
\DeclareMathOperator{\std}{std}

\renewcommand{\epsilon}{\varepsilon}
\newcommand{\bphi}{\phi}

\newcommand{\beq}[1]{\begin{equation}\label{#1}}
\newcommand{\eeq}{\end{equation}}

\newcommand{\RR}{\mathbb{R}}

\newcommand{\derlie}{\mathcal{L}} 
\newcommand{\liederivative}{L}
\newcommand{\ip}[2]{\langle #1,#2 \rangle} 
\newcommand{\lie}[2]{[#1,#2]} 
\newcommand{\courant}[2]{\llbracket#1,#2\rrbracket} 

\newcommand{\db}{\circ} 
\newcommand{\anchor}{\rho} 

\newcommand{\thalf}{\tfrac{1}{2}}
\newcommand{\rond}{\smalcirc} \newcommand{\smalcirc}{\mbox{\tiny{$\circ$}}}

\newcommand{\gendex}[2]{\left\{ #1 \right\}_{#2}}

\newcommand{\cinf}{C^{\infty}} 
\newcommand{\sections}[1]{\Gamma(#1)}
\newcommand{\vf}{\mathfrak{X}} 
\newcommand{\df}{\Omega} 

\newcommand{\ortho}{^{\perp}} 
\newcommand{\ann}{^0} 
\newcommand{\inv}{^{-1}}

\newcommand{\graded}{^{\scriptscriptstyle\bullet}}

\newcommand{\DD}{\mathcal{D}}

\newcommand{\der}{\delta}

\newcommand{\couranta}{(E,\rho,\courant{\cdot}{\cdot},\ip{\cdot}{\cdot})}
\newcommand{\loday}{(E,\rho,\db,\ip{\cdot}{\cdot})}

\newcommand{\kerr}{\ker{\rho}}
\newcommand{\gm}{\Gamma}
\newcommand{\be }{\begin{eqnarray*}}
\newcommand{\eee }{\end{eqnarray*}}

\newcommand{\nah}[1]{H_{\naive}^{#1}(E)}

\newcommand{\ee}{S}
\newcommand{\cala}{\mathcal{A}}
\newcommand{\nD}{\nabla}

\newcommand{\exto}[1]{\xrightarrow{#1}}

\newcommand{\dd}{\breve{d}}
\newcommand{\ii}{\breve{\imath}}

\begin{document}

\title{Modular Classes of Loday Algebroids}

\author{
Mathieu Sti\'enon
\thanks{Research supported by the European Union through the FP6 Marie Curie R.T.N. ENIGMA (Contract number MRTN-CT-2004-5652).} \\
E.T.H.~Zürich \\ 
\texttt{stienon@math.ethz.ch}
\and 
Ping Xu
\thanks{Research partially supported by NSF grant DMS-0605725 \& 
NSA grant H98230-06-1-0047.} \\
Penn State University \\ 
\texttt{ping@math.psu.edu} }

\date{}

\maketitle

\begin{abstract}
We introduce the concept of Loday algebroids, a generalization
of Courant algebroids. We define the \emph{naive cohomology} and modular class 
of a Loday algebroid, and we show that the modular class of 
the double of a Lie bialgebroid vanishes.
For Courant algebroids, we describe the relation between
the naive and standard cohomologies and we conjecture that they are 
isomorphic when the Courant algebroid is transitive.
\end{abstract}

\section{Naive Cohomology}

Given a Courant algebroid $\couranta$, let $\gm(\wedge^k\kerr)$ denote
the space of smooth sections of the (possibly singular) vector bundle
$\wedge^k\kerr$ (i.e. smooth sections $\alpha$ of $\wedge^k E$ such that $\alpha|_m\in\wedge^k\kerr$ for each $m\in M$).
The extension of the pseudo-metric $\ip{\cdot}{\cdot}$ to $\wedge^k E$ 
naturally induces an isomorphism $\Xi:\wedge^k E\to\wedge^k E^*$. 
Since, by definition, $\ip{\DD f}{e}=\thalf\rho(e)f$, the sections of $\wedge^k\kerr$ are characterized as the elements $\epsilon\in\gm(\wedge^k E)$ such that $\ii_{\DD f}\epsilon=0$, $\forall f\in\cinf(M)$. 
Here $\ii_{\DD f}=\Xi\inv\rond i_{\DD f}\rond\Xi$, where $i_{\DD f}:\gm(\wedge^{k+1} E^*)\to\gm(\wedge^k E^*)$ is the usual contraction of exterior forms with the section $\DD f\in\gm(E)$. 
Define an operator  $\dd:\gm(\wedge^k \kerr )\to\gm(\wedge^{k+1}E)$  by 
\begin{multline*} (\dd\alpha )(e_0,\cdots,e_k)=\sum_{i=0}^{k}(-1)^i\rho(e_i)\alpha(e_0,\cdots,\widehat{e_i},\cdots,e_k) \\ +\sum_{i<j}(-1)^{i+j}\alpha(\courant{e_i}{e_j},e_0,\cdots,\widehat{e_i},\cdots,\widehat{e_j},\cdots,e_k) ,\end{multline*} 
for all $\alpha\in\gm(\wedge^k \kerr )$ and $e_0,\dots,e_k\in\gm(E)$.
Here the pairing between $\gm(\wedge^k \kerr )$ and $\gm(\wedge^k E )$
is via the identification $\Xi:\wedge^k E\to\wedge^k E^*$. 
The following Lemma follows from the Courant algebroid properties, in particular the 
relations $\rho(\DD f)=0$ and $\courant{\DD f}{e}+\DD\ip{\DD f}{e}=0$.
\begin{lemma} We have $\dd\gm(\wedge^k\kerr)\subset\gm(\wedge^{k+1}\kerr)$. 
Moreover, $\big(\gm(\wedge\graded\kerr),\dd\big)$ is a cochain complex.
\end{lemma}
The cohomology of this cochain complex is called the naive cohomology of $E$ and is denoted $H_{\naive}\graded(E)$.
\begin{remark}
It is easy to see that a 1-cochain $\theta\in\gm(\kerr)$ is a 1-cocycle if, and only if, 
$\ip{\theta}{\courant{a}{b}}=\rho(a)\ip{\theta}{b}-\rho(b)\ip{\theta}{a}$ for all $a,b\in\gm(E)$, and a 1-coboundary if, and only if, $\theta=\DD f$ for some $f\in\cinf(M)$.
\end{remark}
\begin{remark}
Let $V$ be the $\cinf(M)$-module generated by $\DD\big(\cinf(M)\big)$. 
Since $\ip{\DD f}{a}=\thalf\rho(a)f$, we have $V=\gm\big(\kerr\ortho\big)$ and $\Xi(V)=\gm\big(\kerr\ann\big)=\rho^*\big(\gm(T^*M)\big)$. 
Moreover $V\subset\gm(\kerr)$, for $\rho\rond\DD=0$. 
Therefore, when $\couranta$ is a regular Courant algebroid (i.e. $\rho$ has constant rank), $E/\Xi\inv(\rho^* T^*M)$ is a Lie algebroid and $H_{\naive}\graded(E)$ is the cohomology of this Lie algebroid.
However, in general, $\gm(E)/V$ is only a Lie-Rinehart algebra over $\cinf(M)$. 
One can consider $H_{\naive}\graded(E)$ as its cohomology \cite{Huesmann}.
\end{remark}

\begin{example}
When $E=TM\oplus T^*M$ is an exact Courant algebroid,
$H_{\naive}\graded(E)$ is isomorphic to the  de~Rham cohomology of $M$. 
\end{example}
\begin{example}
If $E$ is a Courant algebroid over a point, i.e. a Lie algebra equipped with a non-degenerate ad-invariant bilinear form, $H_{\naive}\graded(E)$ 
is simply the Lie algebra cohomology.
\end{example}

\section{Relation with standard cohomology}

Courant algebroids can also be obtained as 
derived brackets \cite{YKS}\cite{Roytenberg} using
 degree two super-symplectic manifolds. More precisely,
given a Courant algebroid $\couranta$, $E[1]$ is
a super-Poisson manifold, where the Poisson structure
is induced by the pseudo-metric. There is a minimal
symplectic realization $X\exto{\pi} E[1]$ and a cubic
function $\Theta$ on $X$ such that $\{\Theta,\Theta\}=0$
and, for all $f\in\cinf(M)$ and $e_1,e_2\in\gm(E)$, 
\[ \DD f=\{\Theta,f\} \qquad \text{and} \qquad e_1\db e_2=\{\{\Theta,e_1\},e_2\} ,\] 
where the symbol $\db$ denotes the \emph{asymmetric Dorfman bracket} defined by the relation $a\db b =\courant{a}{b} +\DD\ip{a}{b}$.
Here elements in $\gm (\wedge^k E)$ are viewed
as functions of degree $k$ on $X$ by considering
them as functions on $E[1]$ via the pseudo-metric $\ip{\cdot}{\cdot}$ 
and identifying them with their pull back by $\pi$. Similarly
functions on $M$ are also identified with their pull
back in $X$. By $\cala^k$ we denote the space of
functions on $X$ of degree $k$. Then $(\cala\graded,\{\Theta,\cdot\})$
is a cochain complex. Its cohomology is called the
{\em standard cohomology} by Roytenberg \cite{Roytenberg} and
we shall denote it by $H_{\std}\graded(E)$.

\begin{lemma}
\label{lem:2}
\begin{enumerate}
\item If $c\in\gm(\wedge^k\kerr)$, then $\{\Theta,c\}=\dd c$;
\item If $c\in\gm(\wedge^k E)$ satisfies $\{\Theta,c\}=0$,
then $c\in\gm(\wedge^k\kerr)$ and $\dd c=0$.
\end{enumerate}
\end{lemma}
\begin{proof}
(i) It suffices to prove the case when $k=1$. The general situation
follows from the Leibniz rule.
Now since $\rho(c)=0$, we have $\forall e_1, e_2\in \gm (E)$,
\begin{align*} & \ip{c \db e_2}{e_1}- (\dd c)(e_1, e_2) \\ 
=& \big(-\ip{e_2\db c}{e_1} +2\ip{\DD\ip{c}{e_2}}{e_1}\big)
-\big(\rho(e_1)\ip{c}{e_2}
-\rho(e_2)\ip{c}{e_1} -\ip{c}{\courant{e_1}{e_2}}\big) \\ 
=& \rho(e_2)\ip{c}{e_1}-\ip{e_2\db c}{e_1} 
-\ip{c}{\courant{e_2}{e_1}} \\
=& \rho(e_2)\ip{c}{e_1}-\ip{e_2\db c}{e_1} 
-\ip{c}{e_2\db e_1} \\
=& 0 .\end{align*}
It thus follows that $\{\{\{\Theta,c\},e_2\},e_1\}
-\{\{\dd c,e_2\},e_1\}=0$, which implies
that $\{\Theta,c\}=\dd c$.

(ii) Since $\ii_{\DD f}c=\{\DD f,c\}=\{\{\Theta,f\},c\}
=\{f,\{\Theta,c\}\}=0$ for all $f\in\cinf(M)$, we have $c\in\gm(\wedge^k\kerr)$.
\end{proof}
As a consequence, we have a homomorphism 
$\phi:H_{\naive}\graded(E)\to H_{\std}\graded(E)$.
Lemma~\ref{lem:2} also implies that $\phi$ is an isomorphism
in degrees $0$ and $1$. It is natural to ask when $\phi$ is 
an isomorphism in all degrees. When
$E$ is a Courant algebroid over a point, $\phi$ is
clearly an isomorphism. On the other hand,
when $E$ is the standard Courant algebroid $TM\oplus T^*M$,
both $H_{\naive}\graded(E)$ and $H_{\std}\graded(E)$ are isomorphic to the de~Rham cohomology of $M$. Hence
$\phi$ is also an isomorphism. This leads to the following
\begin{quote}
{\bf Conjecture} When $E$ is a transitive Courant algebroid, 
$\phi$ is an isomorphism.
\end{quote}

\section{Lie derivatives and Loday algebroids}

The Lie derivative of Courant algebroids was introduced
in \cite{Stienon}. Let us recall its definition briefly.
An infinitesimal automorphism of the vector bundle 
$E\exto{\pi}M$ is a vector field on $E$ 
--- a derivation of the algebra $\cinf(E)$ --- 
which preserves the subspaces $\pi^*\cinf(M)$ and $\sections{E}$ 
(whose elements are identified with functions 
linear on the fibers of $\pi$ through the pairing $\ip{\cdot}{\cdot}$).
In other words, it is a covariant differential
operator on $E$, i.e. a pair of differential
operators $\der^0:\cinf(M)\to\cinf(M)$ and
$\der^1:\sections{E}\to\sections{E}$ satisfying 
\[ \der^0(fg) =f\der^0(g)+\der^0(f)g \qquad \text{and} \qquad 
\der^1(fe) =f\der^1(e)+\der^0(f)e, \] 
for any $f,g\in\cinf(M)$ and $e\in\sections{E}$.
It is known \cite{Roytenberg} that
the Lie algebra $\aut(E)$ of infinitesimal automorphisms 
of the Courant algebroid $E$ consists of those 
covariant differential operators $\der=(\der^0,\der^1)$ on $E$
which satisfy the additional properties:
\[ \der^0\ip{e_1}{e_2}=\ip{\der^1 e_1 }{e_2}+\ip{e_1}{\der^1 e_2}  
\qquad \text{and} \qquad 
\der^1\courant{e_1}{e_2}=\courant{\der^1 e_1}{e_2}
+\courant{e_1}{\der^1 e_2} ,\]
for all $e_1,e_2\in\sections{E}$.

For any $e\in\gm(E)$, the pair $\der_e=(\der^0_e,\der_e^1)$ defined by the relations 
$\der^0_e(f)=\anchor(e)f$ and $\der_e^1(x)=e\db x$
is an infinitesimal automorphism of the Courant algebroid $E$, i.e. $\der_e\in\aut(E)$. 
Let us denote the (local) flow generated by the vector field 
on $E$ corresponding to $\der_e$ by $\bphi_t$. 
By abuse of notations, we use the same symbol
$\bphi_t$ (resp. $\bphi_t^*$) to denote its
induced flow on the tensor bundles 
$E_j^i=(\otimes^i E)\otimes(\otimes^j E^*)$ ($i,j\in\gendex{0,1,2,\dots}{}$) 
(resp. the induced action on the spaces of sections of the $E_j^i$'s). 
For any section $\sigma\in\sections{E_j^i}$, define the Lie derivative
$\derlie_z\sigma\in\sections{E_j^i}$ by
$\derlie_z\sigma=\left.\tfrac{d}{d\tau}\bphi_{\tau}^*\sigma\right|_{\tau=0}$.
Thus we have the usual identity:
$\left.\tfrac{d}{d\tau}\bphi_{\tau}^*\sigma\right|_{\tau=t}
=\bphi_t^*(\derlie_z\sigma)$. 
In the following proposition, we give a list of important
properties of this Lie derivative.
\begin{proposition}
\label{pro:lie}
For all $f,g\in\cinf(M)$ and $x,y,z\in\sections{E}$, we have: 
\begin{gather}
\derlie_z f=\rho(z) f \qquad \derlie_z x=z\db x \qquad \derlie_{\DD f}x=0 \qquad \derlie_x\DD f=\DD\derlie_x f
,\label{1} \\ 
\lie{\der}{\derlie_z}=\derlie_{\der^1 z} \qquad\forall\der\in\aut(E) 
,\label{2} \\
\derlie_z(\sigma\otimes\tau)=\derlie_z\sigma\otimes\tau
+\sigma\otimes\derlie_z\tau \qquad\forall\sigma,\tau\in\oplus_{i,j}E^i_j 
,\label{3} \\
\derlie_{\courant{x}{y}}=\lie{\derlie_x}{\derlie_y} 
,\label{4} \\ 
\derlie_z\courant{x}{y}=\courant{\derlie_z x}{y}+\courant{x}{\derlie_z y} 
,\label{5} \\
\derlie_{fx} y=f\;\derlie_x y-\big(\rho(y)f\big)\;x+2\ip{x}{y}\DD f 
,\label{6} \\ 
\derlie_z\ip{x}{y}=\ip{\derlie_z x}{y}+\ip{x}{\derlie_z y} .\label{7} 
\end{gather}
\end{proposition}

\begin{definition}
A Loday algebroid consists of a vector bundle $\pi:E\to M$, a pseudo-metric 
$\ip{\cdot}{\cdot}$ on the fibers of $\pi$, a bundle map $\rho:E\to TM$ and an $\RR$-bilinear operation $\db$ on $\gm(E)$ 
satisfying 
\begin{gather}
e_1\db(e_2\db e_3)=(e_1\db e_2)\db e_3+e_2\db(e_1\db e_3) ,\label{A} \\ 
\rho(e_1\db e_2)=\lie{\rho(e_1)}{\rho(e_2)} ,\label{B} \\
e_1\db(fe_2)=\big(\rho(e_1)f\big)e_2+f(e_1\db e_2) ,\label{C} \\ 
e_1\db e_2+e_2\db e_1=2\DD\ip{e_1}{e_2} ,\label{D} \\ 
\DD f\db e=0 ,\label{E}
\end{gather}
where $\DD:\cinf(M)\to\gm(E)$ is the $\RR$-linear map 
defined by $\ip{\DD f}{e}=\thalf\rho(e)f$.
\end{definition}
\begin{remark}
\begin{enumerate}
\item According to \cite{Uchino}, for Courant algebroids, axioms \eqref{B} and \eqref{C} are redudant. It would be interesting to investigate if it is also the case for Loday algebroids.
\item The Leibniz algebroids studied by several authors \cite{Spain,Gab,Wade} 
are a more general notion. 
\item A Courant algebroid is a Loday algebroid satisfying the additional axiom 
$\rho(e)\ip{e_1}{e_2}=\ip{e\db e_1}{e_2}+\ip{e_1}{e\db e_2}$. 
\end{enumerate}
\end{remark}
\begin{lemma}
If $\loday$ is a Loday algebroid, then
$\rho(\DD f)=0$ and $\courant{\DD f}{e}+\DD\ip{\DD f}{e}=0$, 
for all $f\in\cinf(M)$ and $e\in\gm(E)$. 
Here $\courant{x}{y}=\thalf(x\db y-y\db x)$ as in a Courant algebroid. 
\end{lemma}
\begin{proof}
Applying $\rho$ to both sides of \eqref{D} and making use of \eqref{B}, we get 
$\rho\big(\DD\ip{e_1}{e_2}\big)=0$ for any $e_1,e_2\in\gm(E)$ and thus also 
$\rho\big(\DD\ip{f e_1}{e_2}\big)=0$ for any $f\in\cinf(M)$.
The Leibniz rule $\DD(fg)=g\DD f+f\DD g$ implies that $\rho(\DD f)=0$.
The other relation follows immediately from \eqref{E} and \eqref{D}.
\end{proof}
As a consequence, the definition of the naive cohomology extends from Courant algebroids to Loday algebroids. 

Let $\loday$ be a Loday algebroid. Given a section $z\in\gm(E)$, set 
$\derlie_z f=\rho(z)f$ for $f\in\cinf(M)$ and $\derlie_z x=z\db x$ for $x\in\gm(E)$ and extend $\derlie_z$ to $\gm(\wedge^k E)$ by the Leibniz rule.
\begin{proposition} 
Identities \eqref{4}, \eqref{5} and \eqref{6} still hold for any Loday algebroid.
\end{proposition}
\begin{remark}
It is unknown if the standard cohomology can be 
defined for Loday algebroids. Indeed, it would be
interesting to see if there exists a derived bracket in the sense of Kosmann-Schwarzbach
\cite{YKS} for a Loday algebroid.
\end{remark}

\section{Modular classes}

A Loday algebroid module is a vector bundle $\ee\to M$ endowed
with an $\RR$-linear map 
$ \gm(E)\otimes\gm(\ee)\to\gm(\ee):e\otimes s\mapsto\nD_e s $
satisfying
\begin{align}
& \nD_{\DD f}s=0 && \nD_e(fs)=f\nD_e s+(\rho (e)f)s \label{eq:3} \\
& \nD_{fe}s=f\nD_e s && \nD_{e_1}(\nD_{e_2}s) -\nD_{e_2}(\nD_{e_1}s)=\nD_{\courant{e_1}{e_2}}s \label{eq:4}
\end{align}
for any $f\in\cinf(M)$, $e,e_1,e_2\in\gm(E)$ and $s\in\gm(\ee)$.

Now let $\ee$ be a real line bundle which is a module
of the Loday algebroid $E$. Assume that 
there exists a nowhere zero section $s\in\gm(\ee)$.
The relation $D_e s=\ip{\theta_s}{e}s$ defines a section $\theta_s\in\gm(E)$. 
From $\nD_{\DD f}s=0$, it follows that $\rho(\theta_s)=0$.
And from \eqref{eq:4}, it follows that $\theta_s$ is a naive $1$-cocycle.
Finally, \eqref{eq:3} implies that, for any nowhere vanishing function
$f\in C^{\infty}(M)$, $\theta_{fs}=f \theta_{s}+2\DD(\ln |f|)$. 
Thus the class $[\theta_s]\in \nah{1}$ is independent of the chosen section $s$ and only depends on the module $\ee$. We will denote this class by $\theta_\ee$.
As in \cite{ELW,WHY}, when the line bundle $\ee$ is not trivial, we set 
$\theta_\ee=\frac{1}{2}\theta_{\ee\otimes\ee}$, where $S\otimes S$ is necessarily a trivial real line bundle. We call $\theta_\ee$ the modular class of the module $\ee$.

\begin{theorem}
Given a Loday algebroid $\loday$, $\wedge^{top}E$ is an $E$-module with $\nD=\derlie$.
\end{theorem}
\begin{proof}
It remains to prove that $\derlie_{fe}s=f\derlie_{e}s$ for any $f\in\cinf(M)$ 
and $s\in\gm(\wedge^{top}E)$.
According to Proposition~\ref{pro:lie}, 
for any $f\in\cinf(M)$ and $e,a\in\gm(E)$, we have 
\begin{multline*}\derlie_{fe}a
=f\derlie_{e}a-(\rho(a)f)e+2\ip{e}{a}\DD f
=f\derlie_{e}a-2\ip{\DD f}{a}e+2\ip{e}{a}\DD f \\ 
=\big(f\derlie_{e}-2(e\wedge)\rond\ii_{\DD f}+2(\DD f\wedge)\rond\ii_{e}\big)(a) \end{multline*}
Note that, as differential operators on 
$\gm(\wedge\graded E)$,
$\derlie_{fe}$, $f\derlie_{e}$, $2(e\wedge)\rond\ii_{\DD f}$
and $2(Df\wedge)\rond\ii_{e}$ all are derivations of degree $0$
with respect to the wedge product on $\gm(\wedge\graded E)$.
Since $\derlie_{fe}$ and $f\derlie_{e}-2(e\wedge)\rond\ii_{\DD f}+2(Df\wedge)\rond\ii_{e}$ are equal when acting both on
sections of $E$ and on functions on $M$, they are also equal when extended to 
$\gm(\wedge^* E)$. 
In particular, if $s\in\gm(\wedge^{top}E)$, 
\[ \derlie_{fe} s =f\derlie_{e}s-2(e\wedge)\rond\ii_{\DD f}s+2(\DD f\wedge)\rond\ii_{e}s =f\derlie_{e}s-2\ip{e}{\DD f}s+2\ip{\DD f}{e}s =f\derlie_{e}s .\]
\end{proof}
The modular class $[\theta_{\wedge^{top}E}]\in\nah{1}$ of the $E$-module $\wedge^{top}E$ is called the \emph{modular class of the Loday algebroid} $E$.

\section{Examples}

Let $E=A\oplus A^*$ be the double of a Lie bialgebroid $(A, A^*)$ \cite{LWX}.
In this case, $\DD=\thalf(d+d_*)$ and, for all $X,Y\in\gm(A)$ and $\xi,\eta\in\gm(A^*)$, the bracket on $\gm (E)$ is defined by
\begin{gather*} \courant{X}{\xi}=(-L_{\xi}X+\thalf d_{*}(\xi,X))
+(L_{X}\xi-\thalf d(\xi,X)), \\ 
\courant{X}{Y}= [X,Y], \qquad 
\courant{\xi}{\eta}= [\xi,\eta]. \end{gather*} 

Now $\wedge^{top} E\cong (\wedge^{top}A)\otimes(\wedge^{top}A^*)$
is a trival line bundle. For the sake of simplicity, we assume that there 
exists a nowhere vanishing section $V\in\gm(\wedge^{top}A)$.
Let $\Omega\in\gm(\wedge^{top}A^*)$ be its dual section.
For any $X\in\gm(A)$ and $\xi\in\gm(A^*)$, one has $\derlie_{X}\xi=\courant{X}{\xi}+\DD\ip{X}{\xi}
=-L_\xi X+d_{*}(\xi, X)+L_X\xi =-i_{\xi}d_{*}X+L_{X}\xi$. 
It follows from the Leibniz rule (see Proposition~\ref{pro:lie}) 
that, for any $\sigma\in\gm(\wedge^k A^*)$,
$\derlie_{X}\sigma=L_{X}\sigma+\lambda$
where $\lambda\in\gm\big(A\otimes(\wedge^{k-1}A^*)\big)$.
On the other hand, since $A$ is isotropic with respect to $\ip{\cdot}{\cdot}$, 
we have that $\derlie_{X}\tau=L_{X}\tau$ if $\tau\in\gm(\wedge^k A)$.
It thus follows that
\[ \derlie_{X}(V\wedge \Omega)=(\derlie_{X}V)\wedge\Omega
+V\wedge(\derlie_{X}\Omega)=(L_X V)\wedge\Omega+V\wedge(L_X\Omega)=0 .\]
Similarly we have $\derlie_{\xi}(V\wedge\Omega)=0$, for all
$\xi\in\gm(A^*)$. Thus we have proved
\begin{theorem}
If a Courant algebroid is the double of
a Lie bialgebroid, then its modular class
vanishes.
\end{theorem} 
\begin{example}
If $E=TM\oplus T^*M$ is an exact Courant algebroid, 
the Courant bracket is given by
\[ \courant{X+\xi}{Y+\eta}=\lie{X}{Y}+i_{X\wedge Y}\phi
+\liederivative_{X}\eta-\liederivative_{Y}\xi
+\thalf d\big((\xi,Y)-(\eta,X)\big) ,\] 
for all $X,Y\in\vf(M)$ and $\xi,\eta\in\df^1(M)$. 
Here $\phi$ is a closed 3-form.

Take a nowhere zero $V\in\gm(\wedge^{top}TM)$ 
and its dual $\Omega\in\Omega^{top}(M)$.
One easily sees that $\derlie_{X}V=L_XV+V'$, where $V'\in
\gm (T^*M\otimes\wedge^{top-1}TM)$ and $\derlie_{X} \Omega=
L_X \Omega$. Thus it follows that $\derlie_{X}
(V\wedge \Omega)=L_X (V\wedge \Omega)=0$, $\forall X\in \gm (TM)$. One also 
sees that $\derlie_{\xi} (V\wedge \Omega)=0$, $\forall \xi \in \gm (T^*M)$. 
Therefore the modular class vanishes.
\end{example}

\paragraph{Acknowledgments} 
We thank Zhang-Ju Liu, Alan Weinstein and the referee for useful comments.

\end{document}